\documentclass[a4paper, 12pt]{extarticle}
\usepackage{etex}
\usepackage{verbatim}
\usepackage{enumitem}
\usepackage[latin1]{inputenc}
\usepackage[T1]{fontenc}
\usepackage{lmodern}
\usepackage{textcomp}
\usepackage[italian, english]{babel}
\usepackage{amsmath, amssymb, amscd, amsthm, amsfonts, mathrsfs, imakeidx, tikz}
\makeindex[intoc]

\usepackage[dvips]{geometry}
\DeclareMathOperator{\R}{\mathbb R}    
\DeclareMathOperator{\N}{\mathbb N}    

\newcommand{\extp}{\@ifnextchar^\@extp{\@extp^{\,}}}\def\@extp^#1{\mathop{\bigwedge\nolimits^{\!#1}}}  

\theoremstyle{plain}
\newtheorem{theorem}{Theorem}[section]
\newtheorem{lemma}[theorem]{Lemma}
\newtheorem{proposition}[theorem]{Proposition}
\newtheorem{corollary}[theorem]{Corollary}
\newtheorem{conjecture}[theorem]{Conjecture}
\newtheorem{claim}[theorem]{Claim}
\theoremstyle{remark}
\newtheorem{definition}[theorem]{Definition}
\newtheorem{remark}[theorem]{Remark}

\newtheorem{remarks}[theorem]{Remarks}
\newtheorem{nota}[theorem]{Note}
\newtheorem{example}[theorem]{Example}
\newtheorem{examples}[theorem]{Examples}
\newtheorem{exercise}[theorem]{Exercise}

\newcommand{\bt}{\begin{theorem}}
\newcommand{\et}{\end{theorem}}

\newcommand{\bv}{\begin{vuoto}}
\newcommand{\ev}{\end{vuoto}}

\newcommand{\bl}{\begin{lemma}}
\newcommand{\el}{\end{lemma}}

\newcommand{\bd}{\begin{definition}}
\newcommand{\ed}{\end{definition}}

\newcommand{\beq}{\begin{equation}}
\newcommand{\eeq}{\end{equation}}

\newcommand{\bexa}{\begin{example}}
\newcommand{\eexa}{\end{example}}

\newcommand{\bexas}{\begin{examples}}
\newcommand{\eexas}{\end{examples}}

\newcommand{\bexe}{\begin{exercise}}
\newcommand{\eexe}{\end{exercise}}

\newcommand{\bclaim}{\begin{claim}}
\newcommand{\eclaim}{\end{claim}}

\newcommand{\bprop}{\begin{proposition}}
\newcommand{\eprop}{\end{proposition}}

\newcommand{\bp}{\begin{proof}}
\newcommand{\ep}{\end{proof}}

\newcommand{\bc}{\begin{corollary}}
\newcommand{\ec}{\end{corollary}}

\newcommand{\bq}{\begin{question}}
\newcommand{\eq}{\end{question}}

\newcommand{\bconj}{\begin{conjecture}}
\newcommand{\econj}{\end{conjecture}}

\newcommand{\bproblem}{\begin{problem}}
\newcommand{\eproblem}{\end{problem}}

\newcommand{\bs}{\begin{proof}[Soluzione:]}
\newcommand{\es}{\end{proof}}

\newcommand{\br}{\begin{remark}}
\newcommand{\er}{\end{remark}}

\newcommand{\brs}{\begin{remarks}}
\newcommand{\ers}{\end{remarks}}

\newcommand{\bn}{\begin{nota}}
\newcommand{\en}{\end{nota}}

\newcommand{\st}{\ |\ }

\renewcommand{\epsilon}{\varepsilon}
\renewcommand{\theta}{\vartheta}
\renewcommand{\phi}{\varphi}

\newtheorem{notat}[theorem]{Notation}

\usepackage{hyperref}
\usepackage[all]{xy}
\usepackage{subfig}
\usepackage[labelfont={bf},size={footnotesize}]{caption}
\usepackage[colored]{shadethm}
\usepackage{color, xcolor}
\usepackage{pstricks, colortab}
\usepackage{titlesec}
\usepackage[nottoc]{tocbibind}

\input diagxy
\makeindex

\title {A characterization of distance matrices of weighted hypercube graphs and Petersen graphs}
\author{Elena Rubei, Dario Villanis Ziani}
\date{}

\begin{document}
\maketitle
\subsubsection*{Abstract}
{\small Given a positive-weighted simple connected graph with $m$ vertices, labelled by the numbers $1,\ldots,m$, we can construct an $m \times m$ matrix whose entry $(i,j)$, for any $i,j\in\{1,\dots,m\}$, is the minimal  weight  of a path between $i$ and $j$, where the weight of a path is the sum of the weights of its edges. Such a matrix is called the distance matrix of the weighted graph.
There is wide literature about distance matrices of weighted graphs. 
In this paper we characterize distance matrices of positive-weighted $n$-hypercube graphs. Moreover  we  show that a connected bipartite $n$-regular graph with order $2^n$ is not necessarily the $n$-hypercube graph. Finally
  we give a characterization of distance matrices  of positive-weighted Petersen graphs.}

\section*{Introduction} 
Throughout the paper let $G$ denote a finite simple graph and let $V(G)$ and $E(G)$ 
be respectively its vertex set and its edge set; moreover, we denote by $e(v,w)$
the edge with endpoints $v$ and $w$, if it exists. A graph $G$ endowed with a function $w:E(G)\to \R^+$, where $\R^+$ is the set of the positive real numbers,  is called a \emph{positive-weighted graph}; we denote it  by $\mathcal G=(G,w)$; for any $e \in E(G)$, we call $w(e)$ the \emph{weight} of $e$. \\ Furthermore, for any subgraph $H$, we define $w(H)$ to be the sum of the weights of the edges of $H$. Suppose now that $G$ is connected; the \emph{$k$-weight} of a $k$-subset of vertices $\{v_1,\dots,v_k\}$ is defined to be the minimum among  the weights of the $\mathcal G$-subgraphs whose vertex set contains $v_1,\dots,v_k$. We denote it by
$D_{v_1,\dots,v_k}(\mathcal G)$.
We say that a subgraph $H$ of $G$
\emph{realizes} $D_{v_1,\dots,v_k}(\mathcal G)$ if its
vertex set contains $v_1,\dots,v_k$ and
$w(H)=D_{v_1,\dots,v_k}(\mathcal G)$.
In particular, in the case $k=2$ we can associate to a positive-weighted graph a symmetric matrix whose entries are the distances between two vertices: if we label, in some way,  the vertices  by the numbers  $1,\ldots, m$,  we define the  $(i,j)$-entry of this matrix to be $D_{i,j} ({\mathcal G})$. Obviously, the diagonal entries are zero, while the off-diagonal entries are strictly positive. We can also label only some of the vertices of $G$ 
and consider the matrix whose entries are 
the distances between  labelled vertices.
Such a  matrix is called the  \emph{distance matrix associated to ${\cal G}$ and to the subset of $V(G)$ given by the labelled vertices.} 
The following result characterizes the matrices that are  associated to some positive-weighted graph.
\bt[\bf{Hakimi-Yau}, \cite{HY65}] \label{HY}
A symmetric matrix ${(D_{i,j})}_{i,j\in \{1,\ldots ,m\}}$ with zero diagonal entries and with strictly positive off-diagonal entries is the matrix associated to a positive-weighted graph with vertex set  $\{1,\ldots ,m\}$  if and only if the triangle inequalities hold, that is, if and only if
\[
  D_{i,j}\leq D_{i,k}+D_{k,j}\quad  \forall i,j,k\in \{1,\ldots, m\}.
\]
\et
A square matrix whose  diagonal entries are zero and  the off-diagonal entries are strictly positive is called a \emph{predistance matrix}.
A predistance matrix satisfying the triangle inequalities is called a
\emph{distance  matrix}. By 
Theorem~\ref{HY}, such a matrix is the distance matrix associated to some positive-weighted graph. 
\\
Among the many results on the theory of weighted graphs, we quote also the famous
criterion for a distance matrix to be  the distance matrix associated to a positive-weighted tree and to 
 some subset of its vertex set, see \cite{Bu}, \cite{SimP}, \cite{Za}:
\bt \label{B} {\bf (Buneman-Simoes Pereira-Zaretskii)}
Let $D$ be an $m \times m $ distance matrix.
It is the distance matrix of a positive-weighted tree with vertex set containing $\{1,\ldots , m\}$
if and only if the so called $4$-point condition holds, 
that is,  for all distinct $i,j,k,h  \in \{1,...,m\}$,
the maximum of $$\{D_{i,j} + D_{k,h},D_{i,k} + D_{j,h},D_{i,h} + D_{k,j}
  \}$$ is attained at least twice.
\et
We recall that weighted graphs and their reconstruction from the $2$-weights have applications in several disciplines, such as 
biology,  archaeology, engineering, computer science.  
For instance, biologists represent the evolution of the species by ``phylogenetic trees'', that is, 
positive-weighted trees 
whose vertices represent  species and the weight of an edge is given by how much the DNA sequences (or some segments of the DNA sequences) of the species represented by the vertices of the edge differ; sometimes they also use  graphs that are not trees (because it  may happen that some segments of DNA sequences  are equal for two species with a different evolution history).  
Analogously, also  evolution of manuscripts in archaeology can be represented  by positive-weighted graphs. 
There is a wide literature concerning graphlike dissimilarity families
and treelike dissimilarity families,
in particular about methods to reconstruct weighted trees from their dissimilarity families; these methods are used by biologists to reconstruct phylogenetic trees. 
See for example \cite{DHK12}, \cite{S-S} for overviews on phylogenetic trees. Moreover,
weighted graphs can represent hydraulic webs or railway webs where the weight of an edge is given by the length or the cost 
  (or the difference between the  earnings and the cost) of the line represented by that edge. 
Obviously, weighted graphs can  represent also social or computer networks.
It can be interesting, given a family of real numbers,
to wonder if there exists 
a weighted graph of a particular kind (for instance a  weighted cycle, a weighted caterpillar,    a weighted bipartite graph and so on)  with the given family  as family of  $2$-weights.
A criterion which determines if a family parametrized by the $2$-subsets of a set is the family of $2$-weights of 
a weighted graph of a particular kind
 could help also to detect mistakes in the measurement of 
the $2$-weights in the case we already know
that the weighted graph we deal with is  of particular kind.
For instance, suppose we have 
a hydraulic web (or 
a web to broadcast signals) that can be represented by a hypercube graph 
whose vertices correspond to wells (respectively 
points where we can transmit the signals);
 we define the weight of an edge to be 
the time 
a substance injected in a  vertex of the edge takes to arrive at the other vertex of the edge (respectively 
the time a signal takes to go from a vertex to the other); suppose we can measure the $2$-weights  of any pair of vertices; if they do not satisfy the criterion for a family to be the family of $2$-weights of a weighted hypercube graph, we can deduce that there have been some mistakes in their measurement or that we have supposed something wrong in our model 
(for instance it may be wrong to  assume that the time 
a substance takes to go through  a path is the sum of the periods related to the  edges of the path). Observe that it can be useful to know not only the weighted graph ${\cal G}$ representing the network, but also its subgraph ${\cal G}'$ obtained by removing 
the useless edges, where we say that an edge $e$ is useless if, for any vertices $i,j$ 
of ${\cal G}$, there a path realizing the distance between $i$ and $j$ and not containing $e$: for instance, in a war, if the graph is our enemy's comunication network,
the knowledge of ${\cal G}'$ can be useful to know how to interrupt, at least partially, our enemy's comunications. 
\\
In \cite{BR18} the authors characterized the predistance matrices that are actually distance matrices of some particular graphs, such as paths, caterpillars, cycles, bipartite graphs, complete graphs and planar graphs.\\
We can naturally look for similar results 
for other kinds of graphs.
In this work we give some criterions for a distance matrix to be the distance matrix of  a positive-weighted  $n$-hypercube graph, that is, of a positive-weighted graph whose vertices and edges are respectively the vertices and edges  of the hypercube in $\R^n$. First we consider the case where all the edges are useful, then we consider the general case for $n=3$.
Moreover, in order to show that it was not possible
to deduce easily a characterization of
distance matrices
of $n$-hypercube graphs  from the one  for bipartite graphs by adding a condition equivalent to $n$-regularity,
we  exhibit an example of a connected bipartite $n$-regular graph with order $2^n$ that is not an $n$-hypercube graph. We also exhibit a program to see if a distance matrix is the distance matrix of a hypercube graph.
Finally,  we give a characterization of the distance matrices of  positive-weighted Petersen graphs.

\section{Notation and remarks}

\begin{notat} \label{nota}
$\bullet$   We denote the set of non-negative integers by $\N$  and 
     the set of positive integers by $\N^+$.

$\bullet$     
    For any set $A$, we denote the cardinality of $A$ by $\#A$.

$\bullet$     
 For any matrix $D$, let $D_{(i)}$ be the $i$-th row, let $D^{(j)}$ be the $j$-th column and 
 let ${}^t\! D $ be the transpose of $D$.
  
  $\bullet$ For any graph $G$ and for any $x,y \in V(G)$, let $d(x,y)$ be the minimal number of edges of a path in $G$ with endpoints $x$ and $y$. 
  
\end{notat}

We recall from \cite{BR18} the definitions of indecomposable entry of a distance matrix and of useful edge.

\begin{definition}
  Let $D$ be a distance $m \times m$ matrix for some  $m  \in \N^+$. Let $i,j \in \{1, \dots, m\}$ with $i\neq j$.
  We say that the entry $D_{i,j}$
   is \emph{indecomposable} if and only if $$D_{i,j}<D_{i,k}+D_{k,j}$$ for any  $k \in \{1, \ldots ,m \} \setminus \{i,j\}$. Otherwise we say that $D_{i,j}$ is \emph{decomposable}.
\end{definition}

\begin{definition}
  In a positive-weighted graph ${\mathcal G}$ an edge $e$ is called \emph{useful} if there exists at least one pair of vertices $i$ and $j$ such that all the paths realizing $D_{i,j}(\mathcal G)$ contain the edge $e$.
  Otherwise the edge is called \emph{useless}.
\end{definition}

\br \label{useful}
We recall from \cite{BR18} (Remark 2.3) that, if $D$ is the distance matrix of a positive-weighted graph ${\mathcal G}=(G,w)$, then
$D_{i,j}$
is indecomposable if and only if $E(G)$ contains the edge $e(i,j)$ and $e(i,j)$ is useful; in this case we have that $D_{i,j}({\mathcal G})$ is realized only by the path given only by  the edge
$e(i,j)$ and in particular $w(e(i,j))=D_{i,j}({\mathcal G})$.
\er

\begin{notat} \label{Xkx}
  Let $D$ be a distance $m \times m$ matrix for some $m  \in \N^+ $.
  We denote the set $\{1, \ldots , m\}$ by  $X$ and we fix an element $x $ in  $X$. We can partition the set $X$ as follows:
  \begin{itemize}
    \item let $X_0(x)=\{x\}$;
    \item let $X_1(x)$ be the set of the elements   $y \in X$ such that $D_{x,y}$ is indecomposable; 
    \item let $X_2(x)$ be the set of the elements  $y\in X$ for which the minimum $k$ such that there exist $i_1, \ldots, i_{k-1} \in X$ with  $D_{x,i_1},D_{i_1,i_2}, \ldots,$ $ D_{i_{k-1},y}$ indecomposable, is $2$;
    \item in general, for every $t \in \N^+$, we define $X_t(x)$ to be the set of the elements  $y\in X$ for which the minimum $k$ 
    such that there exist $i_1, \ldots, i_{k-1} \in X$ with  $D_{x,i_1},D_{i_1,i_2}, \ldots, D_{i_{k-1},y}$ indecomposable, is $t$.
  \end{itemize}

Finally, 
for any $A \subset X$ and any $t \in \N$,
we denote   $\bigcup_{a\in A} X_t (a)$ by  $X_t (A)$.
\end{notat}

\section{Distance matrices of weighted hypercube graphs}

In this section we give a characterization of distance matrices of positive-weighted
$n$-hypercube graphs.

\bd Let $n \in \N^+$. The
\emph{$n$-hypercube graph} is the graph $Q_n$ whose vertices and edges  are  respectively the  vertices and the edges of the $n$-hypercube, that is, the graph such that $V(Q_n)={\{0,1\}}^n$ and $E(Q_n)$ is 
\[\left\{e(v,w)\st v,w\in{\{0,1\}}^n  \text{ and } 
\exists !  i \in \{1, \ldots,n\} \text{ such that }v_i\neq w_i\right\}.
\]
\ed
\br Let $Q_n$ be
the $n$-hypercube graph.  Observe that for any  $x,y\in V(Q_n)$,  $$d(x,y)=\# \{i \in \{1, \ldots ,n\}| \; x_i \neq y_i\}.$$ 
 Moreover, we can easily prove that the $3$-hypercube graph is the unique   (connected) bipartite $3$-regular graph. The analogous statement does not hold  for the $n$-hypercube graph with $n\geq 4$; we defer the proof of this fact  to the end of this section.
\er

\br Let ${\cal Q}_n=(Q_n,w)$ be
a positive-weighted $n$-hypercube graph  where each edge is useful and let $X$ denote its vertex set. By Remark \ref{useful}, we
have that $$
  X_k(x)= \{y\in V(Q_n)\ |\ d(x,y)=k\},$$
in fact, by Remark \ref{useful}, the indecomposable $2$-weights correspond to the edges of $Q_n$ and so the  minimum $k$ such that  there exist $i_1, \ldots, i_{k-1} \in X$ with
$D_{x,i_1},D_{i_1,i_2},\dots,D_{i_{k-1},y}$  indecomposable, is equal $d(x,y)$.
\er

We recall now two characterizations of hypercube graphs. Each of them will be useful to give  
a characterization of distance matrices of positive-weighted $n$-hypercube graphs. To state the theorems we need to recall two definitions.

\bd We say that a graph $G$ is a  \emph{$(0,2)$-graph} if, for any $x,y \in V( G) $, the number of the vertices that are adjacent both to $x$  
and to $y$ is either $0$ or $2$.
\ed

\bd A \emph{geodesic} in a graph $G$ is  a path with minimal number of edges, i.e. a geodesic with respect to the distance $d$ defined in Notation \ref{nota}.
\ed

\bt \label{Mulder}  {\bf (Mulder \cite{Muld}, Laborde and Rao Hebbare \cite{Lab-Heb})}
Let $G$ be 
a connected $(0,2)$-graph.
Then it is regular. Let $d$ be its degree. 
We have that  $\#V(G) \leq 2^d$; moreover $G$ is  
a hypercube graph if and only if  $\#V(G) = 2^d$. 
\et

\bt \label{Foldes} {\bf (Foldes \cite{Fol})}
A connected graph $G$ is a hypercube graph if and only if the following two conditions hold:

(i) $G$ is bipartite, 

(ii) for any two vertices $x,y$ of $G$ the number of the geodesics between $x$ and $y$ is $d(x,y)!$. 
\et

We are now ready to give our characterizations of the distance matrices 
of positive-weighted hypercube graphs where all the edegs are useful.

\bt \label{cubici0}
Let $n \in \N^+$  and let $D$ be a  $2^n \times 2^n$ distance matrix.  We denote   the set $\{1, \ldots , 2^n\}$ by $X$.
The matrix $D$ is the distance matrix of a positively weighted $n$-hypercube graph
${\mathcal Q}_n=(Q_n,w)$
in which each edge is useful if and only if
 the following conditions hold:
\begin{enumerate}[label=(\alph*)]
 
  \item  the number of the $\{x,y\} \in {X \choose 2 }$ such that $D_{x,y}$ is  indecomposable   is $ 2^{n-1}n$,
 
 \item if $D_{x,z}$ and $D_{z,y}$  
are indecomposable for some distinct $x,y,z \in X$, 
then there exists exactly one $z' \in X \setminus  \{x,y,z\}$ such that $D_{x,z'}$ and $D_{z',y}$  
are indecomposable.

\end{enumerate}
\et
\bp
$(\Rightarrow)$
By Remark \ref{useful}, the indecomposable $2$-weights correspond to the edges of $Q_n$. So this implication is obvious.

$(\Leftarrow)$
We define $G_n$ to be the graph whose vertex set is $X$ and,  for any $i,j \in X$, we have that  $ E(G_n)$ contains  $e(i,j)$ if and only if $
D_{i,j}$ is indecomposable.

 By (b), the graph 
 $G_n$ is a $(0,2)$-graph. So, by Theorem \ref{Mulder},  it is regular. Let $d$ be its degree. Then $\#E(G) = \frac{(\#X) d}{2}=  \frac{ 2^n d}{2}= 2^{n-1}d
$. By (a), we have that   
 $\#E(G) = 2^{n-1}n $. So we get $d=n$ and then $\# V(G)= 2^d$; therefore, by Theorem \ref{Mulder}, we get that $G_n$ is
 isomorphic to the $n$-hypercube graph.
 
 Let us define $w: E(G_n) \rightarrow \R$ as follows: $w(e(x,y))=D_{x,y}$ for 
every $x,y $ with $D_{x,y}$ indecomposable and let ${\cal G}_n$ be $(G_n, w)$.
We have to prove that $D_{x,y}({\cal G}_n) =
D_{x,y}$ for any $x,y \in X$. If 
$D_{x,y}$ is indecomposable, the statement follows at once from the triangular inequalities. Suppose $D_{x,y}$ is decomposable; then there exist $i_1, \ldots, i_k \in X$ for some $k$ such that $$D_{x,y}=D_{x,i_1}+ D_{i_1, i_2}+ \ldots+ D_{i_k,y},$$ with  
$D_{x,i_1}, D_{i_1, i_2}, \ldots, D_{i_k,y}$ indecomposable. 
 From the equality above and the triangle inequalities we get that  $D_{x,y}$ is equal to $$\min_{{\footnotesize\begin{array}{c}s \in \N\setminus \{0\},\; j_1, \ldots, j_s \in X, \\ D_{x,j_1}, D_{j_1, j_2}, \ldots, D_{j_s,y} \; indecomposable \end{array} }} \{D_{x,j_1}+ D_{j_1, j_2}+ \ldots+ D_{j_s,y}\};$$ the number above is obviously equal to $D_{x,y}({\cal G}_n)$ by the definition of $2$-weights; hence we conclude.
\ep

\bt \label{cubici}
Let $n \in \N^+$  and let $D$ be a  $2^n \times 2^n$ distance matrix.  We denote the set $\{1, \ldots , 2^n\}$ by $X$.
The matrix $D$ is the distance matrix of a positively weighted $n$-hypercube graph
${\mathcal Q}_n=(Q_n,w)$
in which each edge is useful if and only if
 the following conditions hold:
\begin{enumerate}[label=(\alph*)]
  \item  for any $ x, y \in X$ and $k \in \N$ with
         $ y  \in X_k (x) $, we have: $$ \# (X_1 (y) \cap X_{k-1}(x))=k ,$$ 
        \item if $i_1, \ldots, i_k $ are distinct elements of  $X$ and $ D_{i_1,i_2}, \ldots, D_{i_{k-1},i_k} $   
are indecomposable, 
then there do not exist distinct  $j_1, \ldots, j_{k-1} \in X$ such that $D_{j_1,j_2}, \ldots, D_{j_{k-2},j_{k-1}} $   
are indecomposable with $i_1=j_1$, $i_k=j_{k-1}$.
\end{enumerate}
\et
\bp
$(\Rightarrow)$
By Remark \ref{useful}, the indecomposable $2$-weights correspond to the edges of $Q_n$.
So
$X_k(x) $ is the set of the $n$-tuples with entries in $\{0,1\}$ with $n-k$ entries equal to the corresponding entries of $x$  and the others different from the corresponding entries of $x$. So
(a) is obvious.
Also (b) is obvious,
in fact, if two paths in a hypercube graph have the same endpoints, 
 the number of the edges of one is odd if and only if  the number of the edges of the other is odd.

$(\Leftarrow)$
We define $G_n$ to be the graph whose vertex set is $X$ and,  for any $i,j \in X$, we have that $e(i,j) \in E(G_n)$ if and only if $i\in X_1 (j)$.
We want to  show that $G_n$ is isomorphic to the $n$-hypercube graph.

First let us prove  that (b) implies that 
\begin{equation} \label{pp}
X_1(
  X_r (\overline{x})) \cap
  X_s (\overline{x}) = \emptyset \hspace*{1cm} \mbox{\rm for } r+s  
 \; \mbox{\rm even} \end{equation}
  for  any $ \overline{x}$ vertex of $G_n$. Since  $X_1(
  X_r (\overline{x})) \cap
  X_s (\overline{x}) \neq \emptyset$
if and only if  $X_1(
  X_s (\overline{x})) \cap
  X_r (\overline{x}) \neq \emptyset$, we can suppose $r \leq s$. 
  If $r < s$, then $ r \leq s-2$ (since $r+s$ is even) and thus the statement (\ref{pp}) is obvious from our definition of
   $X_s (\overline{x})$. 
  Let $r=s$; suppose there exists $z \in  
X_1( X_r (\overline{x})) \cap  X_s (\overline{x}) =   X_1( X_r (\overline{x})) \cap
  X_r (\overline{x}) $; since $z \in  
X_1( X_r (\overline{x})) $, then there exist distinct $y, j_1, \ldots, j_{r-1}$ such that  $D_{z,y}$,
  $D_{y,j_1}$,  $D_{j_1,j_2}$,...,  $D_{j_{r-1},
  \overline{x}}$ are indecomposable; 
  since $z $ is also in $ X_r (\overline{x})$, we must have that $z$ is different from any of 
  $ j_1, \ldots, j_{r-1}$  (if not, $z$ would be in $X_l (\overline{x})$  for some $l < r$); moreover $z \neq y$ because $y \in X_1(z)$; 
  hence $z,y, j_1, \ldots, j_{r-1}$  are distinct and such that  $D_{z,y}$,
  $D_{y,j_1}$,  $D_{j_1,j_2}$,...,  $D_{j_{r-1},
    \overline{x}}$ are indecomposable; 
   since $z $ is also in $ X_r (\overline{x})$,  we get a contradiction by (b). 

Formula (\ref{pp})
 implies  that $G_n$ is bipartite: 
 for any $\overline{x} \in V(G_n)$, the graph $G_n$ is bipartite on  the sets $$ \cup_{k \; \mbox{\rm\footnotesize even}} X_k(\overline{x}), \;\;\;\;\;\; \cup_{k \; \mbox{\rm\footnotesize odd}} X_k(\overline{x}) ,$$ in fact:
let $ y,y' $ be in the first set; by (\ref{pp}) there is no edge in $G_n$ 
 with endpoints $y,y'$, i.e. 
 $y' \not \in X_1 (y)$; analogously if 
  $ y,y' $ are in the second set. 
  
Now let us prove that, for any $x$ and $y$ vertices of $G_n$,  the number of the geodesics between $x$ and $y$ is $d(x,y)!$.
 Suppose $y \in X_s (x)$, hence $d(x,y)=s$.
 Let $i_1, \ldots, i_{s-2}$ be vertices of $G_n$ such that $D_{x,i_1}$,...., $D_{i_{s-2}, y}$
 are indecomposable, hence such that $x,i_1, \ldots, i_{s-2}, y$ are vertices of a geodesic between $x$ and $y$. By (a) with $k=s$, we have $s$ choices for $i_1$; if we fix $i_1$, then by (a) with $k=s-1$ we have $ 
 s-1 $ choices for $ i_2$ and so on. Hence the number of the geodesics is $s!$.
 
By Foldes' theorem we can conclude that 
$G_n$  is isomorphic to the $n$-hypercube graph. Let us define $w: E(G_n) \rightarrow \R$ as follows: $w(e(x,y))=D_{x,y}$ for 
every $x,y $ with $y \in X_1(x)$.
 Let ${\cal G}_n$ be $(G_n, w)$.
We can prove that $D_{x,y}({\cal G}_n) =
D_{x,y}$ for any $x,y \in X$ as in the final part of the proof of Theorem \ref{cubici0}. 
\ep

\bigskip
We want now to exhibit
a program that allows us to see if a
distance matrix is the distance matrix of a hypercube graph where all the edges are useful. We use the concept of \emph{tower matrix} introduced in \cite{SP}:
given an $m \times m$ matrix $D$, its tower matrix $T$ is a $m^2 \times m$ matrix 
whose rows are indicized by the elements $(i,j)$ with $i,j \in \{1, \ldots, m\}$ ordered by lexicographic order and whose $(i,j)$-row is 
$D_{(i)}+{}^t \! D^{(j)}$. Let  $D$ be a distance matrix; then 
the $l$-entry of  the $(i,j)$-row is the minimal weight of a path between $i$ and $j$ and passing through $l$; observe that $D_{i,j}$ is indecomposable if and only if in the   $(i,j)$-row of $T$ the minimum is reached only twice. Here we see the tower matrix as a tridimensional $m \times m \times m$ matrix $T$ such that $T(j,\cdot, i) = D_{(i)}+{}^t \! D^{(j)}$.
We see if an $m \times m$
distance matrix $D$ is the distance matrix of a hypercube graph where all the edges are useful by checking if $\#X_1(x) = \log_2(m)$ for any $x \in \{1, \ldots,m\}$ and if the graph given by the indecomposable entries of $D$ is a $(0,2)$-graph.
 The program requires $O(m^3)$ elementary operations.
 
{\small
\begin{verbatim}
function H=hypercube(D)

% D is a symmmetric matrix with positive off-diagonal entris and zero diagonal
% entries  and the entries satisfy the triangle inequalities;
% H will be 1 if D is the distance matrix of a hypercube graph where all the
% edges are useful, 0 otherwise

m=size(D,2);
d=log2(m)
H=1;
T=zeros(m,m,m);
IND=zeros(m,m);

for i=1:m
        T(:,:,i)=D+ones(m,1)*D(i,:);
        IND(:,i)=  (sum(  (T(:,:,i)==D(:,i)).')).'== 2*ones(m,1);
% First we calculate the i-part of the tower matrix T seen as 3 dimensional 
% matrix. Then we search for the indecomposable entries of D.
% IND(:,i) is a column of length m whose j entry is 1 when D_ij is
% indecomposable  and it is 0 when D_ij is decomposable with i different 
% from j; the i entry is 0
      if sum (IND(:,i))~=d
      H=0
% in fact if D is the distance matrix of a hypercube graph where all the edges 
% are useful then for every i the number of the indecomposable D_ij 
% must be the logarithm of m in base 2
      return
     end
end

% Now we check if the graph given by the indecomposable D_i,j 
% is a (0,2) graph
for i=1:m
for j=i+1:m                 
                  if ismember(IND(i).'*IND(j),[0.2])=0
                  H=0;
                  return
                  end
end
end

                  \end{verbatim}
}

\bigskip
We want now to state a theorem 
 characterizing distance matrices
of positive-weighted $3$-hypercube graphs
without the assumption that every edge is useful. First we need to recall some notation
and a theorem.

\bd Let $D$ be a distance $m \times m$  matrix. We say that it is \emph{median} if, for any $a,b,c \in \{1, \ldots, m\}$, there exist a unique element $ y  \in \{1, \ldots, m\}$ such that $$D_{i,j}= D_{i,y} + D_{y,j}$$ for any distinct $i,j \in \{a,b,c\}$.
\ed

\bd We denote by $St(n_1, \ldots, n_k)$ 
the starlike tree with the degree of the root equal to $k$ and the number of the edges of the paths between the root and the leaves 
equal to $n_1, \ldots, n_k$.

We denote by $Br(s|n_1, \ldots, n_k| m_1, \ldots, m_h) $ the tree satisfying the following conditions:
\begin{itemize}
\item
there are only two vertices of degree greater than $2$; let us call them  $v_1$ and $v_2$,
\item
the number of the edges of the path $p$  between $v_1$ and $v_2$ is equal to $s$, 
\item
the paths  with endpoint $v_1$ and different from $p$ have   $n_1, \ldots, n_k$ edges,  
\item 
the paths  with endpoint $v_2$ and different from $p$ have   $m_1, \ldots, m_h$ edges.  
\end{itemize}
Obviously any permutation of $n_1, \ldots, n_k$ and any permutation of $m_1, \ldots, m_h$ give the same tree; moreover, $Br(s|n_1, \ldots, n_k| m_1, \ldots, m_h) $ is isomorphic to
$Br(s|m_1, \ldots, m_h| n_1, \ldots, n_k) $. 
\ed

\bd
Let $D$ be  an $m \times m$ distance matrix. A \emph{path $p$ of indecomposable entries} is a sequence $i_1, \ldots, i_k$ in $\{1, \ldots, m\}$ 
such that $D_{i_j, i_{j+1}}$ is indecomposable
for every $j=1,\ldots, k-1$.  We call $k-1$ the \emph{length} of $p$.
\ed

The following theorem, perhaps well-known to experts, can be found for instance in 
 \cite{H-F} and in 
 \cite{B-R5}.

\bt \label{med}
Let $D$ be an $m \times m$ distance matrix. There exists a positive-weighted tree ${\cal T}=(T,w)$, with $V(T)=\{1,\ldots, m\}$, such that $D_{i,j}({\cal T})=D_ {i,j}$ for every $i,j \in \{1,\ldots, m\}$ if and only if the $4$-point condition holds and $D$ is median.
\et

\br \label{tree7}
Any connected subgraph $S$ 
of $Q_3$ such that $V(S) =V(Q_3)$
and with exactly $7$ edges  is a tree (if it contained a cycle, then it would be a $t$-cycle with $t \geq 4$ and $t$ even;
since $\#E(S)=7$, $t$ can be only $4$ or $6$,
 but then we could not have $V(S)=V(Q_3)$).
 Moreover,  we can easily see that a tree with $7$ edges and the degree of any vertex less than or equal to $3$, can be only one of the following: 
\begin{itemize}
\item the path with $7$ edges, 

\item
the starlike trees $St(1,2,4)$,   $ St(1,3,3)$,
$ St(1,1,5)$, $ St(2,2,3)$,

\item
the trees $Br(1| 2,2| 1,1)$, $Br(1| 1,2| 1,2)$, $Br(1| 1,3| 1,1)$,

\item
the tree $Br(2|1,1|1,2)$.
\end{itemize}

We can easily see that they  all are subgraphs of $Q_3$  apart from the following:
  $$ St(1,3,3),\; \;  St(1,1,5), \; \; Br(1| 2,2| 1,1),\;\; Br(2|1,1|1,2).$$ In particular we can deduce that a tree with $7$ edges and the degree of every vertex less than or equal to $3$ is  a subgraph of $Q_3$ if and only if 
  there exists a maximal path with odd length and, in case  there are two vertices of degree $3$, then the path between them 
  has only one edge.
\er

\br \label{seven}
 A connected subgraph $S$ 
of $Q_3$ such that $V(S) = X$ 
 must have at least $7$ edges: let $p$ be a path in $S$  connecting  a vertex with its opposite vertex; the number of edges of $p$ must be odd and  greater than or equal to  $3$; if it is $3$ we need at least $4$ edges to connect $p$ with the elements of $X$ which are not vertices of $p$; 
 if it is $5$ we need at least $2$ edges to connect $p$ with the elements of $X$ which are not vertices of $p$.
\er

\bt
Let  $D$ be a  $8 \times 8$ distance matrix.  Let us denote the set $\{1, \ldots , 8\}$ by $X$.
The matrix $D$ is the distance matrix of a positive-weighted $3$-hypercube graph
${\mathcal Q}_3=(Q_3,w)$ if and only if
 the following conditions hold:
\begin{enumerate}[label=(\alph*)]
  \item  
  $\# X_1(x) \leq 3$ for any $x \in X$;
  
  \item the number $r$ of the $\{x,y\} \in {X \choose 2 }$ such that $D_{x,y}$ is  indecomposable   is in the set $   \{7, \ldots,12\}$;
  
  \item
  \subitem $\bullet$  if $r=7$, the matrix  $D$ satisfies the $4$-points condition, is median,  there exists a maximal path of indecomposable entries with odd length and, if there are two elements in $X$, $z_1 $ and $z_2$,  with $\# X_1(z_1) = \# X_1(z_2) =3 $,  then $D_{z_1, z_2}$ is indecomposable;

  \subitem  $\bullet$ if $r \in \{8, 9, 10,11\}$, then 
  we have: 
  \subsubitem   {\bf --} there exists  $z_1, \ldots,  z_k \in X$ for some $k \in \{\min\{13-r,4\}, \ldots, \min\{24-2r,6\}\}$  such that  
  $\#X_1(x) =3 $ if and only if 
   $x \in X \setminus \{z_1,\ldots, z_k\}$ 
 
  \subsubitem {\bf --} for every $i=1, \ldots,k$ there is  a $(3-\# X_1(z_i))$-subset $S_i$ of $\{z_1, \ldots, z_k\} \setminus \{z_i\}$  such that, if we define 
   $ \tilde{X}_1(z_i) = X_1(z_i) \cup S_i$ and  
    we define $
  \tilde{X}_1 (x)= X_1(x) $ $  \forall x \in X \setminus \{z_1, \ldots, z_k\} $, we have that 
  $$\#(\tilde{X}_1(x) \cap \tilde{X}_1(y) ) 
  \in \{0,2\}$$ for any distinct $x,y \in X$;

  \subitem  $\bullet$ if $r=12$, then  $\#(X_1(x) \cap X_1(y) ) 
  \in \{0,2\}$ for any distinct $x,y \in X$.
  \end{enumerate}
\et

\bp 
Let $S$  be the graph whose vertex set is $X$ and such that $e(i,j)$ is an edge if and only if $D_{i,j}$ is indecomposable, that is $i \in X_1(j) $, and 
 let ${\mathcal S}=(S,w_S)$, where $w_S(e(i,j))=D_{i,j}$ for any $i,j$ with $ i \in X_1(j)$. 

Observe that $S$ is connected: let $x , y \in X$ with $x \neq y$; if $D_{x,y}$ is indecomposable, then there is an edge with endpoints $x$ and $y$; if 
$D_{x,y}$ is decomposable, then it can be written as sum of indecomposable entries of $D$, so  there is a path  with endpoints $x$ and $y$.

First we want to see that
 $S$ is a subgraph of $Q_3$ if and only if conditions (a), (b) and  (c) hold.
 By arguing as in the final part of the proof of Theorem \ref{cubici0}, we get that $D$ is the distance matrix of ${\cal S}$.
Suppose that   $S$ is a subgraph of $Q_3$; 
then, from Remark \ref{seven}, we can deduce  (b); moreover it is obvious that also (a) holds. In the case $r=7$,  we can deduce 
 (c)  from Remark  \ref{tree7} and Theorem \ref{med}; in the cases $r=8,9,10,11,12$, it is obvious that (c) holds.
On the other side, suppose (a), (b) and (c) hold. If  $r=7$, then $S$ is a a subgraph of $Q_3$  by Theorem \ref{med} and Remark \ref{tree7}.
Suppose $r \in \{8,9,10,11\}$; the graph $S$  is the subgraph of the graph $Q$ whose vertex set is $X$ and such that $e(i,j)$ is an edge if and only if  $i \in \tilde{X}_1(j) $; by Theorem \ref{Mulder} we have that $Q$ is isomorphic to $Q_3$. So $S$ is a subgraph of $Q_3$. 
If $r=12$, then $S$ coincides with $Q_3$.

Now we prove the two implications of the statement of the theorem.

From Remark \ref{useful}, we can easily prove the following fact ($\ast$): if 
two edges, $e$ and $e'$,  in  a positive-weighted graph 
${\cal G}$ are useless, then 
$e'$ is useless also in the graph we obtain from ${\cal G}$ by removing $e$.
Suppose that  $D$ is the distance matrix of a positive-weighted $3$-hypercube graph
${\mathcal Q}_3 =(Q_3,w)$;
then, by ($\ast$), it is the distance matrix of 
the subgraph obtained from ${\mathcal Q}_3$ 
by removing  the useless edges, which coincides with ${\mathcal S}$.
So $S$ is a subgraph of $Q_3$ and then (a), (b), (c) hold. 

Conversely, suppose
that (a), (b), (c) hold; we have proved that  this implies that $S$ is a subgraph of $Q_3$. 
Let ${\cal Q}_3 =(Q_3,w)$ be the weighted graph such that $w(e(i,j)) =w_S(e(i,j)) =D_{i,j}$ for any $i,j \in X$ such that $e(i,j) \in S$, i.e. 
$D_{i,j}$ is indecomposable, and such that
\begin{equation} \label{bo}
 w(e(i,j))= 2 \max_{x,y \in X}\{D_{x,y}\}
 \end{equation}
  for any distinct $i,j$ such that $e(i,j) \not\in E(S)$.
 Hence $$D_{i,j} ({\mathcal Q}_3)= D_{i,j} ({\mathcal S}),$$ 
 where the inequality 
 $D_{i,j} ({\mathcal Q}_3) \leq D_{i,j} ({\mathcal S}) $ is obvious, while the inequality $D_{i,j} ({\mathcal Q}_3) \geq  D_{i,j} ({\mathcal S}) $ follows from (\ref{bo}). We have already proved that $D_{i,j} ({\mathcal S}) =D_{i,j}$, therefore we get 
 $D_{i,j} ({\mathcal Q}_3)= D_{i,j} $ for any distinct $i$ and $j$.
\ep

As we have already said, the distance  matrices of positive-weighted bipartite graphs were  characterized in \cite{BR18}.
Obviously an $n$-hypercube graph is a  $n$-regular bipartite graph with order $2^n$.
In order to show that it was not possible
to deduce  a characterization of
distance matrices
of $n$-hypercube graphs  from the one  for bipartite graphs
simply by adding the condition that, for
every vertex $x$, there are exactly $n$ other vertices $y_1, \ldots y_n$ such that
$D_{x, y_i}$ is indecomposable,
we show an example of a connected bipartite $n$-regular graph with order $2^n$ that is not an $n$-hypercube graph. We exhibit here the case $n=4$, being the general case completely analogous.

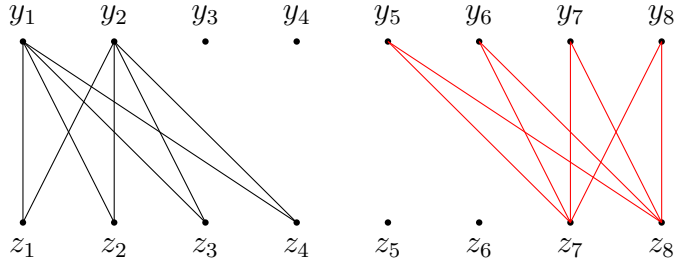
\begin{figure}[h]
  \centering
  {\begin{tikzpicture}[scale=0.60]
      \foreach \x/\xtext in {0/y_1,2/y_2,4/y_3,6/y_4,8/y_5,10/y_6,12/y_7,14/y_8}{
          \fill (\x,0) circle (2pt) node [above=2pt]{$\xtext$};}
      \foreach \x/\xtext in {0/z_1,2/z_2,4/z_3,6/z_4,8/z_5,10/z_6,12/z_7,14/z_8}{
          \fill (\x,-4) circle (2pt) node [below=2pt]{$\xtext$};}
      \foreach \y in {0,2,4,6}{
          \foreach \x in {0,2}
          \draw (\x,0)--(\y,-4);}
      \foreach \y in {12,14}{
          \foreach \x  in {8,10,12,14}
          \draw [red] (\x,0)--(\y,-4);}
    \end{tikzpicture}}
  \caption{Construction of $K_{4,2}$}
  \label{fig:K42}
\end{figure}

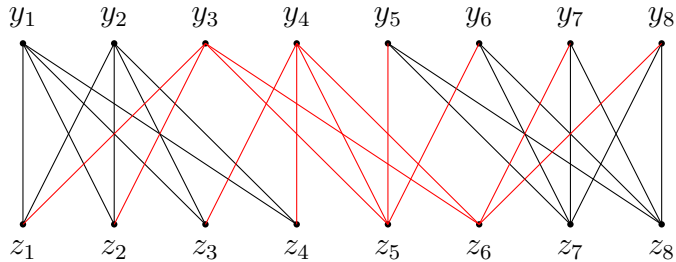
\begin{figure}[h]
  \centering
  {\begin{tikzpicture}[scale=0.60]
      \foreach \x/\xtext in {0/y_1,2/y_2,4/y_3,6/y_4,8/y_5,10/y_6,12/y_7,14/y_8}{
          \fill (\x,0) circle (2pt) node [above=2pt]{$\xtext$};}
      \foreach \x/\xtext in {0/z_1,2/z_2,4/z_3,6/z_4,8/z_5,10/z_6,12/z_7,14/z_8}{
          \fill (\x,-4) circle (2pt) node [below=2pt]{$\xtext$};}
      \foreach \y in {0,2,4,6}{
          \foreach \x in {0,2}
          \draw (\x,0)--(\y,-4);}
      \foreach \y in {12,14}{
          \foreach \x  in {8,10,12,14}
          \draw (\x,0)--(\y,-4);}
      \foreach \y in {0,2,8,10}{
          \draw [red] (4,0)--(\y,-4);}
      \foreach \y in {4,6,8,10}{
          \draw [red] (6,0)--(\y,-4);}
      \foreach \x in {8,10}{
          \draw [red] (\x,0)--(8,-4);}
      \foreach \x in {12,14}{
          \draw [red] (\x,0)--(10,-4);}
    \end{tikzpicture}}
  \caption{Links for the vertices which were still isolated}
  \label{fig:freever}
\end{figure}

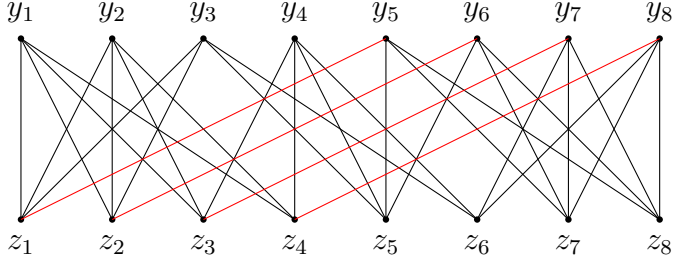
\begin{figure}[h]
  \centering
  {\begin{tikzpicture}[scale=0.60]
      \foreach \x/\xtext in {0/y_1,2/y_2,4/y_3,6/y_4,8/y_5,10/y_6,12/y_7,14/y_8}{
          \fill (\x,0) circle (2pt) node [above=2pt]{$\xtext$};}
      \foreach \x/\xtext in {0/z_1,2/z_2,4/z_3,6/z_4,8/z_5,10/z_6,12/z_7,14/z_8}{
          \fill (\x,-4) circle (2pt) node [below=2pt]{$\xtext$};}
      \foreach \y in {0,2,4,6}{
          \foreach \x in {0,2}
          \draw (\x,0)--(\y,-4);}
      \foreach \y in {12,14}{
          \foreach \x  in {8,10,12,14}
          \draw (\x,0)--(\y,-4);}
      \foreach \y in {0,2,8,10}{
          \draw (4,0)--(\y,-4);}
      \foreach \y in {4,6,8,10}{
          \draw (6,0)--(\y,-4);}
      \foreach \x in {8,10}{
          \draw (\x,0)--(8,-4);}
      \foreach \x in {12,14}{
          \draw (\x,0)--(10,-4);}
      \foreach \x in {8,10,12,14}{
          \draw [red] (\x,0)--(\x-8,-4);}
    \end{tikzpicture}}
  \caption{A connected bipartite $4$-regular graph with order $2^4$ that is not a $4$-hypercube graph.}
  \label{fig:graph}
\end{figure}

To construct our example,
we start by partitioning the set of vertices $X$ (of cardinality equal to $16=2^4$), into two equipotent subsets of cardinality $8$, say $Y=\{y_1,\dots,y_8\}$ and $Z=\{z_1,\dots,z_8\}$.\\
Now we build the complete bipartite graph $K_{2,4}$ with vertex set $$\{y_1,y_2,z_1,z_2,z_3,z_4\},$$ connecting $y_1$ and $y_2$ to each $z_i$
and
  we build the complete bipartite graph $K_{4,2}$ with vertex set $$\{y_5,y_6,y_7,y_8,z_7,z_8\}$$ connecting $z_7$ and $z_8$ to each $y_j$; now $z_1$ and $z_2$ have degree $4$, while each $y_j$ for $j=5, \ldots, 8$ has degree $2$: see Figure \ref{fig:K42}.\\
At this point $y_3$, $y_4$, $z_5$ and $z_6$ are still ``isolated''; since we are building a $4$-regular graph, $z_1$, $z_2$, $z_3$, $z_4$, $y_5$, $y_6$, $y_7$ and $y_8$ can be linked to two more vertices each. So we connect $y_3$ to $z_1$, $z_2$, $z_5$ and $z_6$, then we connect $y_4$ to $z_3$, $z_4$, $z_5$ and $z_6$; moreover we build an edge from $z_5$ to $y_5$ and $y_6$ and from $z_6$ to $y_7$ and $y_8$, as in Figure \ref{fig:freever}.
In this situation, the vertices $y_5$, $y_6$, $y_7$, $y_8$, $z_1$, $z_2$, $z_3$ and $z_4$ have degree $3$, while the others have degree $4$; so we simply connect $z_i$ to $y_{i+4}$ for each $i\in\{1,2,3,4\}$, having a connected bipartite $n$-regular graph with order $2^n$, as desired (Figure \ref{fig:graph}).\\
But this graph is not a $4$-hypercube graph, since there exist at least two different vertices (for example $y_1$ and $y_2$) connected to the same four vertices, while this does not happen in an $4$-hypercube graph. In fact, given two different vertices $u=(u_1,u_2,u_3,u_4)\in\{0,1\}^4$ and $v=(v_1,v_2,v_3,v_4)\in\{0,1\}^4$ of a $4$-hypercube graph, they have at least one different coordinate, say $0=u_1\neq v_1=1$;
we can suppose $u=(0,0,0,0)$; if $u$ and $v$  were adjacent to the same vertices $x,y,z,w$, then each of them would have exactly one coordinate different from $u$, say $0=u_1\neq x_1=1$, $0=u_2\neq y_2=1$, $0=u_3\neq z_3=1$, $0=u_4\neq w_4=1$, so $x=(1,0,0,0)$, $y=(0,1,0,0)$, $z=(0,0,1,0)$, $w=(0,0,0,1)$;
but also $v$ must have the same property,
and since $v\neq u$, we have necessarily $v=(1,1,0,0)$ or $v=(1,0,1,0)$ or $v=(1,0,0,1)$, using the adjacency of $v$ and $x$;
in each case, we see that $v$ and at least one of $y,z$ and $w$   differ in two coordinates, so  they cannot be adjacent
and this is a contradiction.

\section{Distance matrices of weighted  Petersen graphs}

In this section we characterize distance matrices
of positive-weighted Petersen graphs.

\bt \label{petersen}
Let $D$ be a $10 \times 10$ distance matrix.
We denote the set $\{1, \ldots, 10\}$ by  $X$.
The matrix $D$   is the distance matrix of a positively weighted Petersen graph in which each edge is useful if and only if the following conditions hold:
\begin{enumerate}[label=(\alph*)]
  \item
        for any $x \in X$, we have $\#X_1(x)=3$;
  \item for any   $k \in \{3,4\}$ and any  $i_1, \ldots , i_k \in X$ such that  $ D_{i_1, i_2}, D_{i_2, i_3}, \ldots ,$ $ D_{i_{k-1}, i_k}
        $ are all indecomposable, we have that  $D_{i_1, i_k}$ is decomposable;
  \item there exist distinct $v_1,\ldots, v_5  \in X$ such that $D_{v_1, v_2}, \ldots, D_{v_4, v_5}, D_{v_5, v_1}$ are indecomposable and,
        if we denote by $\overline{v_j}$ the unique element in $X_1(v_j)  \setminus \{v_1, \ldots , v_5\}$, we have that $\overline{v_1}, \ldots, \overline{v_5}$ are distinct.
\end{enumerate}
\et
\bp
First observe that, if (a) and (b) hold, then, for any $v_1, \ldots, v_5 $ as in (c), we have that
$\# X_1(v_j) \setminus \{v_1, \ldots ,v_5\}=1 $ for any $j=1, \ldots ,5$: by (b), $D_{v_i, v_j}$
is indecomposable if and only if either  $j=i \pm 1 $ or
$\{i,j\}=\{1,5\}$, so $X_1(v_j) \ni v_i $ if and only if either
$j=i \pm 1 $ or
$\{i,j\}=\{1,5\}$; moreover, by (a), we have that  $\#X_1(v_j) =3$ for any $j=1, \ldots ,5$, so we conclude.

\begin{figure}[h!]
  \begin{center}

    \begin{tikzpicture}
      \draw [thick] (0,1) --(0,2);
      \draw [thick] (-0.951,0.309) --(-1.902 ,0.618);
      \draw [thick] (-0.587,-0.809) --(-1.174 ,-1.618);
      \draw [thick] (0.587,-0.809) --(1.174 ,-1.618);
      \draw [thick] (0.951,0.309) --(1.902 ,0.618);

      \draw [thick] (0,2)--(-1.902 ,0.618);
      \draw [thick] (-1.902 ,0.618)--(-1.174 ,-1.618);
      \draw [thick] (-1.174 ,-1.618)--(1.174 ,-1.618);
      \draw [thick] (1.174 ,-1.618)-- (1.902 ,0.618);
      \draw [thick] (1.902 ,0.618)--(0,2);

      \draw [thick] (0,1) --(-0.587,-0.809);
      \draw [thick] (-0.951,0.309) --(0.587,-0.809) ;
      \draw [thick] (-0.587,-0.809) --(0.951,0.309);
      \draw [thick] (0.587,-0.809) --(0,1);
      \draw [thick] (0.951,0.309) --(-0.951,0.309);

      \node[above] at (0,2) {$v_1$};
      \node[left] at (-1.902 ,0.618) { $v_2$};
      \node[below] at (-1.174 ,-1.618) {$v_3$};
      \node[below] at (1.174 ,-1.618) { $v_4$};
      \node[right] at (1.902 ,0.618) { $v_5$};

      \node[above] at (0.2,1) {$\overline{v_1}$};
      \node[above] at (-0.951,0.309) { $\overline{v_2}$};
      \node[below] at (-0.56,-0.809) {$\overline{v_3}$};
      \node[below] at  (0.587,-0.809) { $\overline{v_4}$};
      \node[above] at (0.951,0.309)  { $\overline{v_5}$};

    \end{tikzpicture}

    \caption{ Petersen graph \label{fig:P}}
  \end{center}
\end{figure}
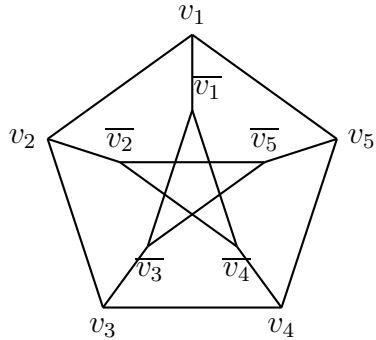

$(\Rightarrow)$
We point out that  all the edges are useful by assumption and an edge $e(i,j)$ is useful
if and only if $D_{i,j}({\mathcal G})$ (which is equal to
$D_{i,j}$ by assumption)   is indecomposable (see Remark \ref{useful}); so the edges correspond to the indecomposable $2$-weights. Hence
statement  (a) follows from the fact that all the vertices of the Petersen graph have degree $3$ and
statement (b) follows from the fact that in the Petersen graph there are not cycles of length $3$ or $4$. Statement  (c) is obvious (take $v_1, \ldots  v_5$ as in Figure \ref{fig:P}).

$(\Leftarrow)$ Let $v_1, \ldots, v_5 , \overline{v_1}, \ldots, \overline{v_5}$ as in (c). Let $G$ be the graph in Figure
\ref{fig:P}  and for any adjacent vertices $i, j $,  define $w(e(i,j))=D_{i,j}$. Let ${\mathcal G} =
  (G,w)$. We want to show that
$D_{x,y}({\mathcal G})=D_{x,y}$ for any $x,y \in X$.

\emph{Case 1}: $D_{x,y}$ is indecomposable.\\
First observe that, by (b), we have that   $D_{v_i, v_j}$ is indecomposable if and only if $\{i,j\} $ is one of the following: $\{1,2\}, \{2,3\}, \{3,4\}, \{4,5\}, \{5,1\}$. So we observe that   $D_{v_i, v_j}$ is indecomposable if and only if
in the graph $G$   we have constructed there is  an edge with endpoints $v_i$ and $v_j$.\\
Observe also that $D_{\overline{v_i}, \overline{ v_j}}$ is indecomposable if and only if in the graph $G$  there is an edge with endpoints $\overline{v_i}$ and $ \overline{ v_j}$; otherwise we would have a
contradiction with (b):
for instance, if $D_{\overline{v_1}, \overline{ v_2}}$ were indecomposable, we would have that
$D_{\overline{v_1},  v_1}$,   $D_{v_1,  v_2}$,
$D_{v_2, \overline{ v_2}}$,
$D_{\overline{v_2}, \overline{ v_1}}$ are   indecomposable (in fact
$D_{\overline{v_1},  v_1}$ is
indecomposable because $\overline{v}_1 \in X_1(v_1) $,
$D_{v_1,  v_2}$  is indecomposable by assumption (c),   and,  finally, $D_{v_2, \overline{ v_2}}$  is indecomposable because $\overline{v}_2 \in X_1(v_2) $)
and this would contradict
assumption (b).\\
Finally
observe that $D_{\overline{v_i}, v_j}$ is indecomposable if and only if in $G$ there is an edge with endpoints $\overline{v_i}$ and $ v_j$
(otherwise we would have again a contradiction with  (b)).\\
Thus we can conclude that $D_{x,y}$ is indecomposable
if and only if in the graph $G$ we have constructed there is an edge with endpoints $x$ and $y$.
In this case we
have that $D_{x,y}
  ({\mathcal G})=
  D_{x,y}$ by the triangle inequalities.

\emph{Case 2}: $D_{x,y}$ is decomposable.   
Then there exist $i_1, \ldots, i_k \in X$ for some $k$ such that $$D_{x,y}=D_{x,i_1}+ D_{i_1, i_2}+ \ldots+ D_{i_k,y}.$$  
 From the equality above and the triangle inequalities we get that  $D_{x,y}$ is equal to $$
 \min_{{\footnotesize\begin{array}{c}s \in \N\setminus \{0\},\; j_1, \ldots, j_s \in X, \\ D_{x,j_1}, D_{j_1, j_2}, \ldots, D_{j_s,y} \; indecomposable \end{array} }} 
  \{D_{x,j_1}+ D_{j_1, j_2}+ \ldots+ D_{j_s,y}\};$$ the second member is  equal to $D_{x,y}({\cal G}_n)$ by the definition of $2$-weights,
by the fact (we have proved before) that $D_{i,j}$ is indecomposable if and only if
in $G$ there is an edge with endpoints $i$
and $j$ and by the fact that $w(e(i,j))=D_{i,j}$ 
for any edge $e(i,j)$.

\ep

\section{Open problems}

We list here some possible open problems.

(1) We could try to generalize the result for
positive-weighted Petersen graphs in Section 3
to positive-weighted Kneser graphs.

(2)
Let $n$ be a natural number with $n \geq 2$ and let $\{m_I\}_{I  \in {\{1,...,n\} \choose 2}}$ and
$\{M_I\}_{I  \in {\{1,...,n\} \choose 2}}$ be
two families of positive real numbers    with $m_I \leq M_I$ for any $I$;
in the paper \cite{Ru}
the author studied when
there exist a positive-weighted graph ${\mathcal G}$
and an  $n$-subset $\{1,..., n\}$ of its vertex set such that
$D_I ({\mathcal G}) \in [m_I, M_I] $ for any $I  \in {\{1,...,n\} \choose 2}$
and the  analogous problem for trees. It would be interesting to study when
there exist a positive-weighted
graph of a particular kind (for instance a hypercube, a cycle, a bipartite graph...)
and an  $n$-subset $\{1,..., n\}$ of its vertex set such that
$D_I ({\mathcal G}) \in [m_I, M_I] $ for any $I $. 

(3) In the last years $k$-weights of weighted graphs for $k \geq 3$ have been investigated, see for instance \cite{BR16}, \cite{HHMS12}, \cite{PS04}.   One could try to characterize
families of $k$-weights of some particular graphs for $k \geq 3$.

\bigskip

{\bf Acknowledgments.}
This work was supported by the National Group for Algebraic and Geometric Structures and their  Applications (GNSAGA-INdAM).

{\small
\phantomsection
}

{\bf Address}: Dipartimento di Matematica e Informatica ``U. Dini'',
viale Morgagni 67/A,
50134  Firenze, Italia

rubei@math.unifi.it, dario.villanisziani@unifi.it
\end{document}